\documentclass{amsart}

\usepackage{ifxetex}
\ifxetex
\usepackage{fontspec}
\else
\usepackage[utf8]{inputenc}
\fi

\usepackage[parfill]{parskip}    
\usepackage{graphbox}
\usepackage{amsmath, amsfonts, amssymb}
\usepackage{url}	
\usepackage{caption}
\usepackage{subcaption}

\usepackage{graphicx}
\usepackage{bbold}

\usepackage{biblatex}
\newcommand\s[1]{s_{#1}}
\newcommand\si[1]{\overline{s_{#1}}}

\title{Weighted growth functions of automatic groups}
\author{Mikael Vejdemo-Johansson}
\date{\today}

\begin{document}

\begin{abstract}
  The growth function is the generating function for sizes of spheres around the identity in Cayley graphs of groups. 
  We present a novel method to calculate growth functions for automatic groups with normal form recognizing automata that recognize a single normal form for each group element, and are at most context free in complexity:
  context free grammars can be translated into algebraic systems of equations, whose solutions represent generating functions of their corresponding non-terminal symbols.
  This approach allows us to seamlessly introduce weightings on the growth function: assign different or even distinct weights to each of the generators in an underlying presentation, such that this weighting is reflected in the growth function.
  We recover known growth functions for small braid groups, and calculate growth functions that weight each generator in an automatic presentation of the braid groups according to their lengths in braid generators.
\end{abstract}

\maketitle

\section{Introduction}

Analytic combinatorics provides tools for enumerating structures as described by formal grammars, producing generating functions.
In this paper we will approach the enumeration of minimal length words representing group elements in finitely presented automatic groups using generating functions generated from the formal grammars associated to the group's automatic structure.

For a group with presentation $G=\langle g | r \rangle$, we define
 \begin{description}
 \item[Cayley graph] the graph with group elements as vertices, and an edge from each vertex $h$ for each generator in $g$, to the vertex $gh$.
 \item[geodesic word] shortest word in the generators and their inverses representing a group element; corresponds to a shortest path in the Cayley graph.
 \item[radius $r$ sphere around the identity] the set of elements whose geodesic words have length $r$. We denote this $S(r)$. 
 \item[growth function] the generating function of the sequence $S(r)$ for $r$ non-negative integers.
 \end{description}

First, in Section~\ref{sec:count-with-gramm} we will introduce the route from a formal grammar to a generating function, and in Section~\ref{sec:braids} we will demonstrate how these methods apply to automatic group, by working with the explicit example of the braid group $B_3$ on three strands.

\section{Counting with grammars}
\label{sec:count-with-gramm}

Chomsky and Sch\"utzenberger proved~\cite{chomsky1963algebraic} that a contextfree language can studied using generating functions.
Their article provides a construction for finding the generating function related to a specific grammar.

Starting with a Backus-Naur form of the grammar, each rewriting rule can be translated into an algebraic equation.
Each terminal symbol is assigned some expression in the variables of the resulting generating function, and each non-terminal symbol is assigned a generating function of its own.
The rewriting assignment is replaced by an equality, each concatenation with a multiplication and each disjunction with an addition. 

For a first and simple example, balanced two-symbol sequences have the grammar
\[
S \to \emptyset \;\; | \;\; \texttt{a}S\texttt{b}
\]

Translating this to an algebraic equation, we would get
\[
S(x,y) = 1 + xy S(x,y)
\]
by weighting each symbol \texttt{a} by $x^1$ and each symbol \texttt{b} by $y^1$. 
The resulting generating function will count the number of strings by the number of \texttt{a} and \texttt{b} symbols in the result, or by evaluating $S(t,t)$ will count by length of the string.

This equation is solved straightforwardly to $S(x,y) = 1/(1-xy) = \sum (xy)^j$, from which we can immediately read that there is exactly one string for each combination of $j$ each of \texttt{a}s and \texttt{b}s.
From $S(t,t) = 1/(1-t^2) = \sum t^{2j}$ follows that there is one unique string for each even length, and no odd-length strings.

Chomsky and Sch\"utzenberger proved that as long as the grammar is at most context-free, the corresponding generating function(s) will be rational functions. 

For anything that can be described by a context-free grammar, this suggests a concrete approach for enumeration: 
\begin{enumerate}
\item Find a Backus-Naur form of a grammar describing your structures
\item Translate the grammar to a system of polynomial equations
\item Use a Gr\"obner basis with an elimination order to solve the system of equations
\item Isolating the Gr\"obner basis elements concentrated to the interesting non-terminal symbol and the terminal variables, solve for a rational form of the generating function
\end{enumerate}

\section{Braids and Automatic Groups}
\label{sec:braids}

\begin{figure}
\includegraphics[]{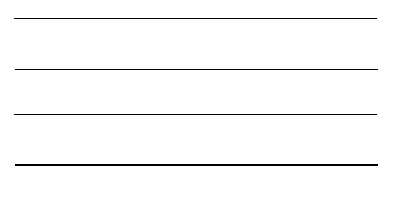}
\includegraphics[]{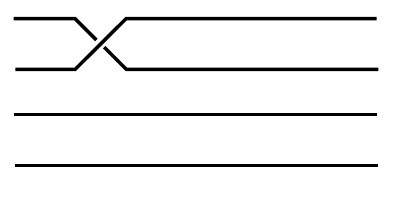} \\[1cm]
\includegraphics[]{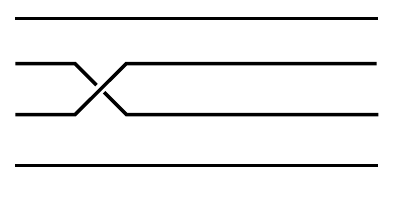}
\includegraphics[]{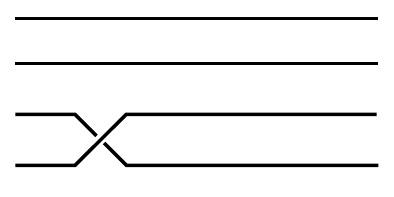} \\[1cm]
\caption{Generators of the Braid group $B_4$}\label{fig:braids-gen}
\end{figure}

\begin{figure}
\includegraphics[align=c]{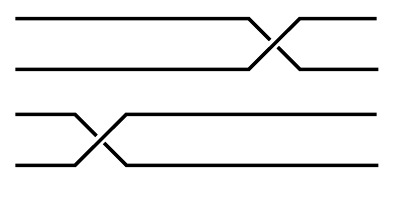} = 
\includegraphics[align=c]{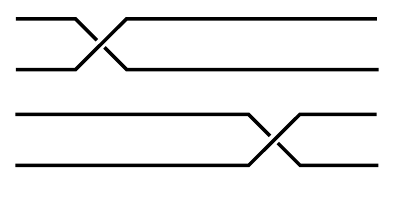} \\[1cm]
\includegraphics[align=c]{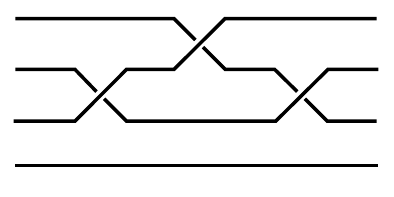}  = 
\includegraphics[align=c]{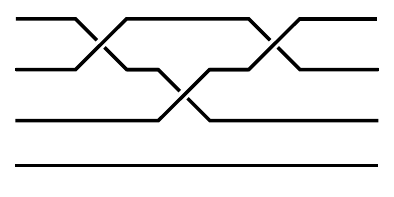}
\caption{Relations of the Braid group $B_4$}\label{fig:braids-rel}
\end{figure}

Braid groups are usually introduced with a finite presentations in terms of elementary braids: for $k$ strands, the braid group $B_k$ has generators $\sigma_j$ for $1\leq j < k$, where $\sigma_j$ crosses strand $j$ over strand $j+1$. We give an illustration for $B_4$ in Figure~\ref{fig:braids-gen}. By inspecting the effects of Reidemeister moves, and of manipulations of separated areas of the 3-sphere, we can derive the finite presentation
\[
B_k = \langle \sigma_1, \dots, \sigma_k | 
  \sigma_i\sigma_j = \sigma_j\sigma_i; 
  \sigma_i\sigma_{i+1}\sigma_i = \sigma_{i+}\sigma_i\sigma_{i+1} \rangle
\]

where $|i-j| > 1$.  Figure~\ref{fig:braids-rel} shows these relations in $B_4$, the smallest braid group where all relations are applicable.

One proof that the word problem is solvable for braid groups was described in~\cite{epstein1992word}, demonstrating that automatic groups solve the word problem, and that braid groups are (bi)automatic.
An automatic group, here, is a finitely presented group coupled with several automata: one to detect whether a given string in the generators is the normal form of a group element, and one to detect right products of a normal form by a generator.

Braid groups form an example of \emph{biautomatic groups}: there are grammars both for recognizing right products and left products.
The part that really interests us here, though, is the normal form recognizer.
With a grammar for normal forms, algebraic equations to compute generating functions for group sizes can be computed.
These generating functions are also studied extensively for finitely presented groups: they are called \emph{growth series}.
For braid groups, we even know grammars that pick out \emph{exactly} one normal form for each group element, such that this normal form is \emph{geodesic}: has the shortest possible expression in some specific set of generators.

\citeauthor{charney1995geodesic} in \cite{charney1995geodesic} gives a grammar for the Braid group $B_3$ with the following transition rules, with each state terminal.

\begin{align*}
B_3 &\to  e | 
 \s1 v_2 | \s2 v_3 | \s1\s2 v_3 | \s2\s1 v_2 | \s1\s2\s1 v_1 | 
 \si1 v_5 | \si2 v_6 | \si1\si2 v_6 | \si2\si1 v_5 | \si1\si2\si1 v_4 \\
v_1 &\to  \s1\s2\s1 v_1 \\
v_2 &\to  \s1 v_2 | \s1\s2 v_3 | \si2 v_6 | \si2\si1 v_5 | \s1\s2\s1 v_1\\
v_3 &\to  \s2 v_3 | \s2\s1 v_2 | \si1 v_5 | \si1\si2 v_6 | \s1\s2\s1 v_1\\
v_4 &\to  \si1\si2\si1 v_4 | \si1\si2 v_6 | \si2 v_6 | \si2\si1 v_5 | \si1 v_5 \\
v_5 &\to  \si1 v_5  | \si1\si2 v_6 \\
v_6 &\to  \si2 v_6 | \si2\si1 v_5 \\
\end{align*}

The construction of this automaton generalizes to all braid groups, with exponential growth in the number of rules in the grammar.
This improves on previous constructions that needed factorial growth in the number of rules.

From this grammar we can produce a system of algebraic equations that counts each generator in Charney's presentation equally
\begin{align*}
B_3 &=(1+ t\cdot v_2 + t\cdot v_3 + t\cdot v_3 + t\cdot v_2 + t\cdot v_1 + 
 t\cdot v_5 + t\cdot v_6 + t\cdot v_6 + t\cdot v_5 + t\cdot v_4 ) \\
v_1 &=(1+ t\cdot v_1 ) \\
v_2 &=(1+ t\cdot v_2 + t\cdot v_3 + t\cdot v_6 + t\cdot v_5 + t\cdot v_1) \\
v_3 &=(1+ t\cdot v_3 + t\cdot v_2 + t\cdot v_5 + t\cdot v_6 + t\cdot v_1) \\
v_4 &=(1+ t\cdot v_4 + t\cdot v_6 + t\cdot v_6 + t\cdot v_5 + t\cdot v_5 ) \\
v_5 &=(1+ t\cdot v_5  + t\cdot v_6 ) \\
v_6 &=(1+ t\cdot v_6 + t\cdot v_5 ) 
\end{align*}

Solving this system for $B_3(t)$ with term order to eliminate all the $v$s, using your favorite computer algebra system recovers a Gr\"obner basis:
\begin{align*}
&4B_3(t){t}^{3}-8B_3(t){t}^{2}-4{t}^{3}+5B_3(t)t-8{t}^{2}-B_3(t)+5t+1 \\
&-2B_3(t){t}^{2}+3B_3(t)t+2{t}^{2}-B_3(t)+5t+{v_6} \\
&-2B_3(t){t}^{2}+3B_3(t)t+2{t}^{2}-B_3(t)+5t+{v_5} \\
&-20B_3(t){t}^{2}+28B_3(t)t+20{t}^{2}-9B_3(t)+52t+2{v_4}+7 \\
&20B_3(t){t}^{2}-26B_3(t)t-20{t}^{2}+5B_3(t)-54t+6{v_3}-11 \\
&20B_3(t){t}^{2}-26B_3(t)t-20{t}^{2}+5B_3(t)-54t+6{v_2}-11 \\
&4B_3(t){t}^{2}-4B_3(t)t-4{t}^{2}+B_3(t)-12t+6{v_1}-7
\end{align*}

The first of these terms completely avoids all the $v$s, and is the one generator of the elimination ideal. 
This produces a functional equation for $B_3(t)$:
\[
4B_3(t){t}^{3}-8B_3(t){t}^{2}-4{t}^{3}+5B_3(t)t-8{t}^{2}-B_3(t)+5t+1 = 0
\]
which we can rewrite to 
\[
B_3(t)(4t^3-8t^2+5t-1) = 4t^3+8t^2-5t-1
\]
from which follows
\begin{multline*}
B_3(t) = \frac{4t^3+8t^2-5t-1}{4t^3-8t^2+5t-1} = 1 + \frac{2t(8t-5)}{(t-1)(2t-1)^2} \\
1 + 10 t + 34 t^2 + 90 t^3 + 218 t^4 + 506 t^5 + 1146 t^6 + 2554 t^7 + 5626 t^8 + 12282 t^9 + 26618 t^{10} + 57338 t^{11} + \\
\quad 122874 t^{12} + 262138 t^{13} + 557050 t^{14} + 1179642 t^{15} + 2490362 t^{16} + 5242874 t^{17} + 11010042 t^{18} + O(t^{19})
\end{multline*}

This recovers the growth function for $B_3$ as computed by Charney~\cite{charney1995geodesic}.

The method of going through Gr\"obner basis computations, however, is more flexible than Charney's linear algebra approach.
Since we can choose weights at will, we can -- for instance -- compute the growth series of the automatic presentation, as weighted by the number of elementary braid generators used for each word.
Doing this still retains a strong focus on the automatic presentation, and as we will see no longer calculates geodesic (ie shortest) words for the presentation with elementary braid generators.

To achieve this, we weight each term when translating to a system of equations not by the number of automatic generators involved, but by the number of elementary braid generators in each term, producing the system of equations
\begin{align*}
B_3 &=(1+
 t\cdot v_2 + t\cdot v_3 + t^2\cdot v_3 + t^2\cdot v_2 + t^3\cdot v_1 + 
 t\cdot v_5 + t\cdot v_6 + t^2\cdot v_6 + t^2\cdot v_5 + t^3\cdot v_4 ) \\
v_1 &=(1+ t^3\cdot v_1 ) \\
v_2 &=(1+ t\cdot v_2 + t^2\cdot v_3 + t\cdot v_6 + t^2\cdot v_5 + t^3\cdot v_1) \\
v_3 &=(1+ t\cdot v_3 + t^2\cdot v_2 + t\cdot v_5 + t^2\cdot v_6 + t^3\cdot v_1) \\
v_4 &=(1+ t^3\cdot v_4 + t^2\cdot v_6 + t\cdot v_6 + t^2\cdot v_5 + t\cdot v_5 ) \\
v_5 &=(1+ t\cdot v_5  + t^2\cdot v_6 ) \\
v_6 &=(1+ t\cdot v_6 + t^2\cdot v_5 ) 
\end{align*}

Calculating, again, an eliminating Gr\"obner basis produces
\begin{multline*}
B_3(t){t}^{5}+B_3(t){t}^{4}-{t}^{5}-3B_3(t){t}^{3}-{t}^{4}-B_3(t){t}^{2}-{t}^{3}+3B_3(t)t-
{t}^{2}-B_3(t)+t+1 \\
\end{multline*}\begin{multline*}
-4{t}^{4}{B_3(t)}^{3}-8{B_3(t)}^{3}{t}^{3}-2{B_3(t)}^{2}{t}^{4}+8
{B_3(t)}^{3}{t}^{2}-4{B_3(t)}^{2}{t}^{3}+3B_3(t){t}^{4}+12{B_3(t)}^{3}t+ \\
 20{B_3(t)}^{2}{t}^{2}+6B_3(t){t}^{3}+3{t}^{4}-8{B_3(t)}^{3}+30{B_3(t)}^{2}t+8{B_3(t)}^{2}{v_6}+18B_3(t){t}^{2}+6{t}^{3}- \\
 12{B_3(t)}^{2}+27B_3(t)t+12B_3(t){v_6}+6{t}^{2}-6B_3(t)+9t+6{v_6} \\
-2{B_3(t)}^{2}{t}^{4}-4{B_3(t)}^{2}{t}^{3}+4{B_3(t)}^{2}{t}^{2}+B_3(t){t}^{3}+2{t}^{4}+6{B_3(t)}^{2}t+8B_3(t){t}^{2}+3{t}^{3}- \\
 4{B_3(t)}^{2}+10B_3(t)t+4B_3(t){v_6}+4{t}^{2}+2t{v_6}-3B_3(t)+4t+2{v_6}+1 \\
\end{multline*}\begin{multline*}
-2{B_3(t)}^{2}{t}^{4}-4{B_3(t)}^{2}{t}^{3}+3B_3(t){t}^{4}+4{B_3(t)}^{2}{
t}^{2}+5B_3(t){t}^{3}-{t}^{4}+6{B_3(t)}^{2}t+3B_3(t){t}^{2} \\
 -{t}^{3}-4{B_3(t)}^{2}+6B_3(t)t+4B_3(t){v_6}-3{t}^{2}+2{{v_6}}^{2}-B_3(t)-2t+2{v_6}-3 \\
\end{multline*}\begin{multline*}
-B_3(t){t}^{4}-2B_3(t){t}^{3}+{t}^{4}+2B_3(t){t}^{2}+2{t}^{3}+3B_3(t)t+2{t}^{2}-2
B_3(t)+3t+{v_5}+{v_6} \\
\end{multline*}\begin{multline*}
-4B_3(t){t}^{4}-7B_3(t){t}^{3}+4{t}^{4}+7B_3(t){t}
^{2}+7{t}^{3}+9B_3(t)t+9{t}^{2}-6B_3(t)+11t+2{v_4}+4 \\
\end{multline*}\begin{multline*}
56{B_3(t)}^{5}{t}^{4}+88{B_3(t)}^{5}{t}^{3}+24{B_3(t)}^{4}{t}^{4}-80{B_3(t)}^{5}{t}^{2}+40{B_3(t)
}^{4}{t}^{3}-40{t}^{4}{B_3(t)}^{3}-88{B_3(t)}^{5}t- \\
 248{B_3(t)}^{4}{t}^{2}-60{B_3(t)}^{3}{t}^{3}+47{B_3(t)}^{2}{t}^{4}+8{B_3(t)}^{5}-280{B_3(t)}^{4}t+96{B_3(t)}^{4}{v_3}- \\
 256{B_3(t)}^{3}{t}^{2}+109{B_3(t)}^{2}{t}^{3}+14B_3(t){t}^{4}-88{B_3(t)}^{4}-296{B_3(t)}^{3}t+144{B_3(t)}^{3}{v_3}-\\
 282{B_3(t)}^{2}{t}^{2}+25B_3(t){t}^{3}
-101{t}^{4}-148{B_3(t)}^{3}-387{B_3(t)}^{2}t+80{B_3(t)}^{2}{v_3}-385B_3(t){t}
^{2}- \\
 202{t}^{3}+79{B_3(t)}^{2}-570B_3(t)t+12B_3(t){v_3}-168B_3(t){v_6}-205
{t}^{2}+181B_3(t)-303t+6{v_3}-196{v_6}-6 \\
\end{multline*}\begin{multline*}
-868{B_3(t)}^{4}{t}^{
4}-1364{B_3(t)}^{4}{t}^{3}-694{t}^{4}{B_3(t)}^{3}+1240{B_3(t)}^{4}{t}^{2}-1126
{B_3(t)}^{3}{t}^{3}+2011{B_3(t)}^{2}{t}^{4}+ \\
 1364{B_3(t)}^{4}t+4304{B_3(t)}^{3}{t}^
{2}+3545{B_3(t)}^{2}{t}^{3}+1376B_3(t){t}^{4}-124{B_3(t)}^{4}+4846{B_3(t)}^{3}t- \\ 
1488{B_3(t)}^{3}{v_3}+2620{B_3(t)}^{2}{t}^{2}+2595B_3(t){t}^{3}-1825{t}^{4
}+1318{B_3(t)}^{3}+2321{B_3(t)}^{2}t- \\
 2784{B_3(t)}^{2}{v_3}-3915B_3(t){t}^{2}-
3650{t}^{3}+4935{B_3(t)}^{2}-6772B_3(t)t-1216B_3(t){v_3}- \\
 2064B_3(t){v_6}-
3625{t}^{2}+364t{v_3}+4941B_3(t)-5475t-50{v_3}-3336{\it 
v6}-132 \\
\end{multline*}\begin{multline*}
252{B_3(t)}^{4}{t}^{4}+396{B_3(t)}^{4}{t}^{3}-4{t}^{4}{B_3(t)}^{3}-360
{B_3(t)}^{4}{t}^{2}+4{B_3(t)}^{3}{t}^{3}-1077{B_3(t)}^{2}{t}^{4}- \\
 396{B_3(t)}^{4}t-
956{B_3(t)}^{3}{t}^{2}-2039{B_3(t)}^{2}{t}^{3}+170B_3(t){t}^{4}+36{B_3(t)}^{4}-
1084{B_3(t)}^{3}t+ \\
 432{B_3(t)}^{3}{v_3}+1030{B_3(t)}^{2}{t}^{2}+321B_3(t){t}^{3}
+659{t}^{4}-412{B_3(t)}^{3}+1745{B_3(t)}^{2}t+ \\
 456{B_3(t)}^{2}{v_3}+2939B_3(t)
{t}^{2}+1318{t}^{3}-2149{B_3(t)}^{2}+4426B_3(t)t+36B_3(t){v_3}+ \\
 1656B_3(t){
v_6}+1299{t}^{2}-1411B_3(t)+1977t+91{v_2}-53{v_3}+1356
{v_6}-38 \\
\end{multline*}\begin{multline*}
504{B_3(t)}^{4}{t}^{4}+792{B_3(t)}^{4}{t}^{3}-8{t}^{4}{B_3(t)}^{3}-
720{B_3(t)}^{4}{t}^{2}+8{B_3(t)}^{3}{t}^{3}-2154{B_3(t)}^{2}{t}^{4}- \\
 792{B_3(t)}^{4
}t-1912{B_3(t)}^{3}{t}^{2}-4078{B_3(t)}^{2}{t}^{3}+158B_3(t){t}^{4}+72{B_3(t)}^{4}
-2168{B_3(t)}^{3}t+ \\
 864{B_3(t)}^{3}{v_3}+2060{B_3(t)}^{2}{t}^{2}+369B_3(t){t}^{3
}+1500{t}^{4}-824{B_3(t)}^{3}+3490{B_3(t)}^{2}t+ \\
 912{B_3(t)}^{2}{v_3}+6151
B_3(t){t}^{2}+2909{t}^{3}-4298{B_3(t)}^{2}+9125B_3(t)t+72B_3(t){v_3}+ \\
 3312B_3(t){
v_6}+3053{t}^{2}-2822B_3(t)+4409t+182{v_1}-288{v_3}+2712
{v_6}+106
\end{multline*}

The first of these terms is the elimination order projection, producing the functional equation
\[
B_3(t){t}^{5}+B_3(t){t}^{4}-{t}^{5}-3B_3(t){t}^{3}-{t}^{4}-B_3(t){t}^{2}-{t}^{3}+3B_3(t)t-
{t}^{2}-B_3(t)+t+1 = 0
\]
Which we can solve for $B_3(t)$, producing
\[
B_3(t) = \frac{t^5+t^4+t^3+t^2-t-1}{t^5+t^4-3t^3-t^2+3t-1} = 
1 + 4t + 10t^2 + 22t^3 + 44t^4 + 84t^5 + O(t^6)
\]

Comparing this to a hand-enumeration of small braids produces 12 braids using two elementary generators, whereas this enumeration predicts 10. 
The reason for this discrepancy is precisely the fact that geodesic here is measured in terms not of elementary generators but in terms of automatic generators. 
Hence, while $\sigma_1^{-1}\sigma_2$ and $\sigma_2^{-1}\sigma_1$ are both length-2 words in the elementary generating set, they have minimal representatives in the automatic presentation as 
\[
\sigma_1^{-1}\sigma_2 = \sigma_2\sigma_1\sigma_2^{-1}\sigma_1^{-1} = \sigma_1^{-1}D\sigma_2
\qquad\text{and}\qquad
\sigma_2^{-1}\sigma_1 = \sigma_1\sigma_2\sigma_1^{-1}\sigma_2^{-1} = \sigma_2^{-1}D\sigma_1
\]
where $D=\sigma_1\sigma_2\sigma_1\sigma_2^{-1}\sigma_1^{-1}\sigma_2^{-1}$, and hence shows up as length 3 instead. 

\section{Conclusion}
\label{sec:conclusion}

The methods from analytical combinatorics producing generating functions directly from contextfree grammars are directly applicable to the problem of computing growth functions for automatic groups. 
They can be weighted, which provides some insight into how the automatic group geodesic words relate to their presentation in a different choice of generators -- however, for, for instance, braid groups, the automatic presentations tend to sort generators moving the elementary generators to the front and their inverses to the end of a word, which may not produce a geodesic in the simpler presentation.

It is unclear how to get closer to a growth function for the elementary presentation of a braid group.

\printbibliography

\end{document}